\documentclass[10pt]{article}
\usepackage{amssymb,amsmath,graphicx,acronym,amsfonts}

\newtheorem{corollary}{Corollary}[section]
\newtheorem{theorem}{Theorem}[section]
\newtheorem{lemma}{Lemma}[section]
\newtheorem{definition}{Definition}[section]
\newtheorem{proposition}{Proposition}[section]
\newtheorem{example}{Example}[section]
\newtheorem{assum}{Assumption}[section]
\newtheorem{algo}{Algorithm}[section]
\newtheorem{Remark}{Remark}[section]

\def\bc{\begin{corl}}
\def\bc{\end{corl}}
\def\ba{\begin{algo}}
\def\ea{\end{algo}}
\def\br{\begin{Remark}}
\def\er{\end{Remark}}
\def\bs{\begin{assum}}
\def\es{\end{assum}}
\def\bt{\begin{theorem}}
\def\et{\end{theorem}\vskip 3pt}
\def\bl{\begin{lemma}}
\def\el{\end{lemma}}
\def\ep{\end{proposition}}
\def\bp{\begin{proposition}}
\def\qed{\hfill{$\Box$}\vskip 5pt}
\def\be{\begin{example}}
\def\ee{\end{example}}
\def\bd{\begin{definition}}
\def\ed{\end{definition}}
\def\bc{\begin{corollary}}
\def\ec{\end{corollary}}
\def\proof{\noindent\it Proof. \hspace{1mm}\rm}
\input epsf
\begin{document}
\title{\bf Positive Definiteness and Semi-Definiteness of Even Order Symmetric Cauchy Tensors}
\author{Haibin Chen\thanks{Department of Applied Mathematics, The Hong Kong Polytechnic University, Hung Hom,
Kowloon, Hong Kong. Email: chenhaibin508@163.com.}, \quad
Liqun Qi
\thanks{Department of Applied Mathematics, The Hong Kong Polytechnic University, Hung Hom,
Kowloon, Hong Kong. Email: maqilq@polyu.edu.hk. This author's work
was supported by the Hong Kong Research Grant Council (Grant No.
PolyU 502510, 502111, 501212 and 501913).} }

\date{\today}
\maketitle

\begin{abstract} Motivated by symmetric Cauchy matrices, we define symmetric Cauchy tensors and their generating vectors
in this paper.  Hilbert tensors are symmetric Cauchy tensors. An
even order symmetric Cauchy tensor is positive semi-definite if and
only if its generating vector is positive. An even order symmetric
Cauchy tensor is positive definite if and only if its generating
vector has positive and mutually distinct entries.    This extends
Fiedler's result for symmetric Cauchy matrices to symmetric Cauchy
tensors.    Then, it is proven that the positive semi-definiteness
character of an even order symmetric Cauchy tensor can be
equivalently checked by the monotone increasing property of a
homogeneous polynomial related to the Cauchy tensor. The homogeneous
polynomial is strictly monotone increasing in the nonnegative
orthant of the Euclidean space when the even order symmetric Cauchy
tensor is positive definite. Furthermore, we prove that the Hadamard
product of two positive semi-definite (positive definite
respectively) symmetric Cauchy tensors is a positive semi-definite
(positive definite respectively) tensor, which can be generalized to
the Hadamard product of finitely many positive semi-definite
(positive definite respectively) symmetric Cauchy tensors.  At last,
bounds of the largest H-eigenvalue of a positive semi-definite
 symmetric Cauchy tensor are given and several spectral properties on Z-eigenvalues of odd order symmetric Cauchy tensors are shown. Further questions on Cauchy tensors are raised.
\medskip

\vskip 6pt


\end{abstract}

\section{Introduction}
Let $\mathbb{R}^n$ be the $n$ dimensional real Euclidean space and the set consisting of all
natural numbers be denoted by $N$. Suppose $m$ and $n$ are positive natural numbers and denote $[n]=\{1,2,\cdots,n\}$.
A Cauchy matrix (maybe not square) is an $m\times n$ structure matrix assigned to $m+n$ parameters
$x_1,x_2,\cdots, x_m,y_1,\cdots, y_n$ as follows: \cite{Polya}
\begin{equation}\label{e11} C=\left[\frac{1}{x_i+y_j}\right],~i\in [m],j\in [n].\end{equation}
The Cauchy matrix has been studied and applied in algorithm designing \cite{T93,G95,Go94}.
When $x_i=y_i$ in (\ref{e11}), it is a real symmetric Cauchy matrix.
Stimulated by the notion of symmetric Cauchy matrices, we give the following definition.

\bd\label{def11} Let vector $c=(c_1,c_2,\cdots,c_n)\in \mathbb{R}^n$. Suppose that a real tensor $\mathcal{C}=(c_{i_1i_2\cdots i_m})$ is defined by
$$c_{i_1i_2\cdots i_m}=\frac{1}{c_{i_1}+c_{i_2}+\cdots+c_{i_m}},\quad j\in [m],~i_j \in [n].$$
Then, we say that $\mathcal{C}$ is an order $m$ dimension $n$ symmetric Cauchy tensor and the vector $c\in \mathbb{R}^n$ is called the generating vector of $\mathcal{C}$.
\ed

We should point out that, in Definition \ref{def11}, for any $m$ elements $c_{i_1}, c_{i_2}, \cdots, c_{i_m}$ in generating vector $c$, it satisfies
$$c_{i_1}+c_{i_2}+\cdots+c_{i_m}\neq 0,$$
which implies that $c_i\neq 0$, $i\in [n]$.

By Definition \ref{def11}, a dimension $n\times n$ real symmetric Cauchy matrix is an order 2 dimension $n$ real symmetric Cauchy tensor.
It is easy to check that every principal subtensors of a symmetric Cauchy tensor is a symmetric Cauchy tensor with a generating vector being a subvector of the generating vector of the original symmetric Cauchy tensor. In this paper, we always consider $m$-th order $n$ dimensional real symmetric Cauchy tensors. Hence, it can be called Cauchy tensors for simplicity.

Apparently, Cauchy tensors are a class of structured tensors. In recent years, a lot of research papers on structured tensors are presented \cite{Chen14,Ding13,Ding14,HH14,Qi14,qi14,QXX, Song14,song14,YY14,Zhang12}, which studied M-tensors,
circulant tensors, completely positive tensors, Hankel tensors, Hilbert tensors, P tensors and B tensors. These papers not only established results on spectral theory and
positive semi-definiteness property of structured tensors, but also gave some
important applications of structured tensors in stochastic process and data fitting \cite{Chen14,Ding14}. Actually, Cauchy tensors have close
relationships with Hankel tensors and Hilbert tensors.

Suppose Cauchy tensor $\mathcal{C}$ and its generating vector $c$
are defined as in Definition \ref{def11}.  If
$$c_{i_1}+c_{i_2}+\cdots+c_{i_m} \equiv
c_{j_1}+c_{j_2}+\cdots+c_{j_m}$$ whenever
$$i_1 + i_2+\cdots+i_m  =
j_1+j_2+\cdots+j_m,$$ then Cauchy tensor $\mathcal{C}$ is a Hankel
tensor in the sense of \cite{Qi14}. In general, a symmetric Cauchy
tensor is not a Hankel tensor.  If entries of $c$ are defined such
that
$$c_i=i-1+\frac{1}{m},~i\in [n],$$
then Cauchy tensor $\mathcal{C}$ is a Hilbert tensor according to \cite{Song14}.

In this paper, we are interested in the positive semi-definiteness conditions and positive definiteness conditions for even order Cauchy tensors,
in addition, several spectral properties of positive semi-definite Cauchy tensors are given. In the next section, we prove that
an even order Cauchy tensor is positive semi-definite if and only if its generating vector is positive. An even order Cauchy tensor is positive definite
if and only if its generating vector has positive and mutually distinct entries.    This extends Fiedler's result for symmetric
Cauchy matrices in 2010  \cite{Fied10} to symmetric Cauchy tensors.
After that, it is shown that an even order Cauchy tensor is positive semi-definite if and only if there is a monotone increasing homogeneous polynomial
defined in the nonnegative orthant of $\mathbb{R}^n$. And the homogeneous polynomial is strictly monotone increasing when the even order Cauchy tensor is
positive definite. Later in Section 2, we show that the Hadamard product of two positive semi-definite (positive definite respectively)
Cauchy tensors is a positive semi-definite (positive definite respectively) tensor, which can be generalized to the Hadamard product of
finitely many positive semi-definite (positive definite respectively) Cauchy tensors.

In Section 3, several spectral inequalities are presented on the largest H-eigenvalue and the smallest H-eigenvalue of a Cauchy tensor. When a Cauchy tensor
is positive semi-definite, bounds of its largest H-eigenvalue are given. Then, for an odd order Cauchy tensor, we prove that the corresponding
Z-eigenvector is nonnegative if the Z-eigenvalue is positive and the Z-eigenvector is non-positive if the Z-eigenvalue is negative. Furthermore,
if the generating vector of an odd order Cauchy tensor is positive, the Cauchy tensor has nonzero Z-eigenvalues. We conclude this paper with some final
remarks in Section 4.    Some questions are also presented, which can be considered in the future.

By the end of the introduction, we add some comments on the notation that will be used in the sequel.  Vectors are denoted by
italic lowercase letters i.e. $x,~ y,\cdots$, and tensors are written as calligraphic capitals such as
$\mathcal{A}, \mathcal{T}, \cdots.$ Denote $e_i\in \mathbb{R}^n$ as the $i$-th unit coordinate vector for $i\in [n]$,
 and ${\bf0}$ as the zero vector in $\mathbb{R}^n$. Suppose $x=(x_1,x_2,\cdots,x_n),~y=(y_1,y_2,\cdots,y_n)$. Then $x\geq y$ ($x\leq y$)
means $x_i\geq y_i$ ($x_i\leq y_i$) for all $i\in [n]$. If both $\mathcal{A}=(a_{i_1 \cdots i_m})_{1\leq i_j\leq n}$ and $\mathcal{B}=(b_{i_1 \cdots i_m})_{1\leq i_j\leq n}$, $j=1,\cdots,m$,
are tensors, then $\mathcal{A}\geq \mathcal{B}$ ($\mathcal{A}\leq \mathcal{B}$) means $a_{i_1 \cdots i_m} \geq b_{i_1 \cdots i_m}$ ($a_{i_1 \cdots i_m} \leq b_{i_1 \cdots i_m}$)
for all $i_1,\cdots,i_m \in [n]$.

\setcounter{equation}{0}
\section{Positive Semi-definite Cauchy Tensors}

An order $m$ dimension $n$ tensor $\mathcal{A}=(a_{i_1i_2\cdots i_m})$ is called positive semi-definite if for any vector $x\in \mathbb{R}^n$, it satisfies
$$\mathcal{A}x^m=\sum_{i_1,i_2,\cdots,i_m\in [n]}a_{i_1i_2\cdots i_m}x_{i_1}x_{i_2}\cdots x_{i_n}\geq 0.$$
$\mathcal{A}$ is called positive definite if $\mathcal{A}x^m>0$ for all nonzero vector $x\in \mathbb{R}^n$. Similarly, tensor $\mathcal{A}$
is negative semi-definite (negative definite) if $\mathcal{A}x^m\leq 0$ ($\mathcal{A}x^m<0$ for all nonzero vector $x\in \mathbb{R}^n$). It is obvious that minus positive semi-definite
tensors (minus positive definite tensors respectively) are negative semi-definite tensors(negative definite tensors respectively).
Clearly, there is no odd order nonzero positive semi-definite tensors.

In this section, we will give some sufficient and necessary conditions for even order Cauchy tensors to be
positive semi-definite or positive definite. Some conditions are extended naturally
from the Cauchy matrix case.

\bt\label{thm21} Assume a Cauchy tensor $\mathcal{C}$ is of even order. Let $c\in \mathbb{R}^n$ be the generating vector of $\mathcal{C}$. Then Cauchy tensor
$\mathcal{C}$ is positive semi-definite if and only if $c>0$.
\et

\proof For necessity, suppose that an even order Cauchy tensor $\mathcal{C}$ is positive semi-definite.   It is easy to check that all composites of generating vector $c$ are positive since
$$\mathcal{C}e_i^m=\frac{1}{mc_i}\geq0,~i\in [n]$$
where $e_i$ is the $i$-th coordinate vector of $\mathbb{R}^n$. So, $c_i>0$ for all $i\in [n]$, which means $c>0$.

On the other hand, assume that $c>0$.    For any $x\in \mathbb{R}^n$, it holds that
$$\begin{array}{rl} \mathcal{C}x^m=& \sum_{i_1,\cdots,i_m \in [n]}c_{i_1i_2\cdots i_m}x_{i_1}x_{i_2} \cdots x_{i_m}\\
=&\sum_{i_1,\cdots,i_m \in [n]} \frac{x_{i_1}x_{i_2} \cdots x_{i_m}}{c_{i_1}+c_{i_2}+\cdots+c_{i_m}} \\
=&\sum_{i_1,\cdots,i_m \in [n]}\displaystyle \int^1_0 t^{c_{i_1}+c_{i_2}+\cdots+c_{i_m}-1}x_{i_1}x_{i_2} \cdots x_{i_m}dt\\
=&\displaystyle \int^1_0\left(\sum_{i\in [n]}t^{c_i-\frac{1}{m}}x_i\right)^m dt\\
\geq &0.
\end{array}$$
Here the last inequality follows that $m$ is even. By the arbitrariness of $x$, we know that Cauchy tensor $\mathcal{C}$ is positive semi-definite and the
desired result holds. \qed

\bc\label{corol21} Assume that even order Cauchy tensor $\mathcal{C}$ and its generating vector $c\in \mathbb{R}^n$ be defined as in Theorem \ref{thm21}. Then Cauchy tensor
$\mathcal{C}$ is negative semi-definite if and only if $c<0$.
\ec

\bc\label{corol22} Assume that even order Cauchy tensor $\mathcal{C}$ and its generating vector $c\in \mathbb{R}^n$ be defined as in Theorem \ref{thm21}. Then Cauchy tensor
$\mathcal{C}$ is not positive semi-definite if and only if there exist at least one negative element in $c$.
\ec

From the results about H-eigenvalues and Z-eigenvalues in \cite{Qi05}, we have the following result.

\bc\label{corol23} Assume that even order Cauchy tensor $\mathcal{C}$ and its generating vector $c\in \mathbb{R}^n$ be defined as in Theorem \ref{thm21}. If $c>0$, then
all the H-eigenvalues and Z-eigenvalues of Cauchy tensor $\mathcal{C}$ are nonnegative.
\ec

\bt\label{thm22} Assume even order Cauchy tensor $\mathcal{C}$ has generating vector $c=(c_1,c_2,\cdots,c_n)\in \mathbb{R}^n$. Suppose $c_1,c_2,\cdots,c_n$
are positive and mutually distinct, then Cauchy tensor $\mathcal{C}$ is positive definite.
\et
\proof For the sake of simplicity, without loss of generality, assume that
$$0<c_1<c_2<\cdots<c_n.$$
since $c>0$ and $c_1,c_2,\cdots,c_n$ are mutually distinct. From Theorem \ref{thm21}, we know that Cauchy tensor $\mathcal{C}$ is positive semi-definite.

We prove by contradiction that Cauchy tensor $\mathcal{C}$ is positive definite when the conditions of this theorem hold.
Assume there exists a nonzero vector $x\in \mathbb{R}^n$ such that
$$\mathcal{C}x^m=0.$$
By the proof of Theorem \ref{thm21}, one has
$$\displaystyle \int^1_0\left(\sum_{i\in [n]}t^{c_i-\frac{1}{m}}x_i\right)^m dt=0,$$
which means
$$\sum_{i\in [n]}t^{c_i-\frac{1}{m}}x_i\equiv 0,~t\in [0, 1].$$
Thus
$$x_1+t^{c_2-c_1}x_2+\cdots+t^{c_n-c_1}x_n\equiv 0,~t\in (0, 1].$$
By continuity and the fact that $c_1,c_2,\cdots,c_n$ are mutually distinct, it holds that
$$x_1=0$$
and
$$x_2+t^{c_3-c_2}x_3+\cdots+t^{c_n-c_2}x_n\equiv 0,~t\in (0, 1].$$
Repeat the process above.     We obtain
$$x_1=x_2=\cdots=x_n=0,$$
which is a contradiction with $x\neq0$. So, for all nonzero vectors $x\in \mathbb{R}^n$, it holds $\mathcal{C}x^m>0$ and
$\mathcal{C}$ is positive definite. \qed

From this theorem, we easily have the following corollary, which was first proved in \cite{Song14}.

\bc\label{corol234} An even order Hilbert tensor is positive definite.
\ec

From Theorem A of \cite{Fied10}, we know that a symmetric Cauchy matrix
$$C=\left[\frac{1}{c_i+c_j}\right]$$
is positive definite if and only if all the $c_i$'s are positive and mutually distinct. In fact, the theorem below shows that conditions in Theorem \ref{thm22}
is also a sufficient and necessary condition, which is a natural extension of Theorem A of \cite{Fied10}, by Fielder.

\bt\label{thm23} Let even order Cauchy tensor $\mathcal{C}$ and its generating vector $c$ are defined as in Theorem \ref{thm22}. Then, Cauchy tensor
$\mathcal{C}$ is positive definite if and only if the elements of generating vector are positive and mutually distinct.
\et
\proof By Theorem \ref{thm22}, we only need to prove the ``only if'' part of this theorem. Suppose that  Cauchy tensor
$\mathcal{C}$ is positive definite.  Firstly, by
Theorem \ref{thm21}, we know that
$$c_i>0,~i\in [n].$$
We now prove by contradiction that $c_i$'s are mutually distinct. Suppose that two elements of $c$ are equal. Without loss of generality,
assume $c_1=c_2=a>0$. Let $x\in \mathbb{R}^n$ be a vector with elements $x_1=1, x_2=-1$ and $x_i=0$ for the others. Then, one has
$$\begin{array}{rl} \mathcal{C}x^m=& \sum_{i_1,\cdots,i_m \in [n]}c_{i_1i_2\cdots i_m}x_{i_1}x_{i_2} \cdots x_{i_m}\\
=&\sum_{i_1,\cdots,i_m \in [n]} \frac{x_{i_1}x_{i_2} \cdots x_{i_m}}{c_{i_1}+c_{i_2}+\cdots+c_{i_m}} \\
=&\frac{1}{ma}\sum_{i_1,\cdots,i_m \in [2]}x_{i_1}x_{i_2} \cdots x_{i_m}\\
=&\frac{1}{ma} \left[(-1)^m+m(-1)^{m-1}+\frac{m!}{2!(m-2)!}(-1)^{m-2}+\cdots+(-1)^{m-m}\right]\\
=&\frac{1}{ma}\left[1+(-1)\right]^m\\
=&0,
\end{array}$$
where we get a contradiction with the assumption that Cauchy tensor $\mathcal{C}$ is positive definite.
Thus, elements of generating vector $c$ are mutually distinct and the desired result follows.
\qed

We denote the homogeneous polynomial $\mathcal{C}x^m$ as
$$f(x)=\mathcal{C}x^m=\sum_{i_1,\cdots,i_m \in [n]}c_{i_1i_2\cdots i_m}x_{i_1}x_{i_2} \cdots x_{i_m}.$$
For all $x,y \in X\subseteq\mathbb{R}^n$, if $f(x)\geq f(y)$ when $x\geq y$($x\leq y$), we say that $f(x)$ is monotone increasing (monotone decreasing respectively) in $X$. If $f(x)>f(y)$ when $x\geq y$, $x\neq y$($x\leq y$, $x\neq y$), we say that $f(x)$ is strictly monotone increasing (strict monotone decreasing respectively) in $X$.

The following conclusion means that the positive semi-definite property of a Cauchy tensor is equivalent to the monotonicity of a homogeneous polynomial
respected to the Cauchy tensor in $\mathbb{R}^n_+$.
\bt\label{thm24} Let $\mathcal{C}$ be an even order Cauchy tensor with generating vector $c$. Then,
$\mathcal{C}$ is positive semi-definite if and only if $f(x)$ is monotone increasing in $\mathbb{R}^n_+$.
\et
\proof For sufficiency, let $x=e_i, y=0$ and $x\geq y$.   Then we have
$$\frac{1}{mc_i}=\mathcal{C}x^m=f(x)\geq f(y)=\mathcal{C}y^m=0,$$
which implies that $c_i>0$ for $i\in [n]$. By Theorem \ref{thm21}, it holds that Cauchy tensor $\mathcal{C}$
is positive semi-definite.

For necessary conditions, suppose $x,y\in \mathbb{R}^n_+$ and $x\geq y$. Then, we know that
$$\begin{array}{rl}
f(x)-f(y)=&\mathcal{C}x^m-\mathcal{C}y^m\\
=&\sum_{i_1,\cdots,i_m \in [n]}c_{i_1i_2\cdots i_m}(x_{i_1}x_{i_2} \cdots x_{i_m}-y_{i_1}y_{i_2} \cdots y_{i_m})\\
=&\sum_{i_1,\cdots,i_m \in [n]}\frac{x_{i_1}x_{i_2} \cdots x_{i_m}-y_{i_1}y_{i_2} \cdots y_{i_m}}{c_{i_1}+c_{i_2}+\cdots+c_{i_m}}\\
\geq &0.
\end{array}$$
Here, the last inequality follows that $x\geq y$ and the fact that $c_i>0,$ for $i\in [n]$, which means that $f(x)$
is monotone increasing in $\mathbb{R}^n_+$ and the desired result holds. \qed

\bl\label{lema21}Let $\mathcal{C}$ be an even order Cauchy tensor with generating vector $c$. Suppose $\mathcal{C}$
is positive definite, then the homogeneous polynomial $f(x)$ is strictly monotone increasing in $\mathbb{R}^n_+$
\el
\proof From the condition that $\mathcal{C}$ is positive definite, by Theorem \ref{thm23}, we have
$$c_i>0,~i\in [n],$$
where scalars $c_i,~i\in [n]$ are entries of generating vector $c$. For any
$x,y\in \mathbb{R}^n_+$ satisfying that $x\geq y$ and $x\neq y$, there exists index $i\in [n]$ such that
$$x_i>y_i\geq0.$$
Then, it holds that
$$\begin{array}{rl}
f(x)-f(y)=&\mathcal{C}x^m-\mathcal{C}y^m\\
=&\sum_{i_1,\cdots,i_m \in [n],(i_1,i_2,\cdots,i_m)\neq(i,i,\cdots,i)}c_{i_1i_2\cdots i_m}(x_{i_1}x_{i_2} \cdots x_{i_m}-y_{i_1}y_{i_2} \cdots y_{i_m})\\
&+c_{ii\cdots i}(x_i^{m}-y_i^m)\\
=&\sum_{i_1,\cdots,i_m \in [n],(i_1,i_2,\cdots,i_m)\neq(i,i,\cdots,i)}\frac{x_{i_1}x_{i_2} \cdots x_{i_m}-y_{i_1}y_{i_2} \cdots y_{i_m}}{c_{i_1}+c_{i_2}+\cdots+c_{i_m}}\\
&+\frac{1}{mc_i}(x_i^{m}-y_i^m)\\
>&0,
\end{array}$$
which implies that the homogeneous polynomial $f(x)$ is strictly monotone increasing in $\mathcal{R}^n_+$.\qed

Now, we give an example to show that the strictly monotone increasing property for the polynomial $f(x)$ is only a necessary condition for
the positive definiteness property of Cauchy tensor $\mathcal{C}$ but not a sufficient condition.

\begin{example}\label{exam31} Let Cauchy tensor $\mathcal{C}=(c_{i_1i_2i_3i_4})$ with dimension 3 and it has generating vector $c=(1,1,1)$. Then,
$$c_{i_1i_2i_3i_4}=\frac{1}{4},~i_1,i_2,i_3,i_4\in [3]$$
and the homogeneous polynomial
$$f(x)=\mathcal{C}x^4=\frac{1}{4}\sum_{i_1,i_2,i_3,i_4\in [3]}x_{i_1}x_{i_2}x_{i_3}x_{i_4}.$$
By direct computation, we know that $f(x)$ is strictly monotone increasing in $\mathbb{R}^3_+$. From Theorem \ref{thm23}, Cauchy tensor $\mathcal{C}$
is not positive definite. \qed
\end{example}

Let $r_i$ denote the sum of the $i$-th row elements of Cauchy tensor $\mathcal{C}$, which can be written such that
$$r_i=\sum_{i_2,\cdots,i_m \in [n]}\frac{1}{c_i+c_{i_2}+\cdots+c_{i_m}},~i\in [n],$$
where $c=(c_1,\cdots,c_n)$ is the generating vector of Cauchy tensor $\mathcal{C}$. Suppose
$$R=\max_{1\leq i\leq n}r_i,~r=\min_{1\leq i\leq n}r_i.$$ If Cauchy tensor $\mathcal{C}$ is positive semi-definite, by Theorem \ref{thm21},
it is easy to check that
$$R=\sum_{i_2,\cdots,i_m \in [n]}\frac{1}{\underline{a}+c_{i_2}+\cdots+c_{i_m}},~r=\sum_{i_2,\cdots,i_m \in [n]}\frac{1}{\bar{a}+c_{i_2}+\cdots+c_{i_m}},$$
where $\underline{a}=\min_{1\leq i \leq n}c_i,~\bar{a}=\max_{1\leq i\leq n}c_i$.

Now, before giving the next conclusion, we give the definition of eigenvalues of tensors and the definition of irreducible tensors, which will be used in the sequel.

The definition of eigenvalue-eigenvector pairs of real symmetric tensors comes from \cite{Qi05}.
\bd\label{def21} Let $\mathbb{C}$ be the complex field. A pair $(\lambda, x)\in \mathbb{C}\times \mathbb{C}^n\setminus \{0\}$ is called an
eigenvalue-eigenvector pair of a real symmetric tensor $\mathcal{T}$ with order $m$ dimension $n$, if they satisfy
\begin{equation}\label{e21}
\mathcal{T}x^{m-1}=\lambda x^{[m-1]},
\end{equation}
where $\mathcal{T}x^{m-1}=\left(\sum_{i_2,\cdots,i_m=1}^n t_{ii_2\cdots i_m}x_{i_2}\cdots x_{i_m} \right)_{1\leq i\leq n}$ and $x^{[m-1]}=(x_i^{m-1})_{1\leq i\leq n}$ are dimension $n$ vectors.
\ed
In Definition \ref{def21}, if $\lambda\in \mathbb{R}$ and the corresponding eigenvector $x\in \mathbb{R}^n$, we call $\lambda,~x$ H-eigenvalue and
H-eigenvector respectively. Let $\rho(\mathcal{C})$ denote the spectral radius of Cauchy tensor $\mathcal{C}$. Then $\rho(\mathcal{C})>0$ and $\rho(\mathcal{C})$ is an H-eigenvalue of $\mathcal{C}$ when $\mathcal{C}$ is positive semi-definite

The following definition is consistent
with \cite{Chang08} and \cite{Qi13} respectively.
\bd\label{def22} For a tensor $\mathcal{T}$ with order m and dimension n. We say that $\mathcal{T}$ is reducible if there is a nonempty
proper index subset $I\subset [n]$ such that
$$t_{i_1i_2\cdots i_m}=0,~\forall ~i_1\in I,~\forall~i_2,i_3,\cdots,i_m \notin I.$$
Otherwise we say that $\mathcal{T}$ is irreducible.
\ed

\bt\label{thm25} Let even order Cauchy tensor $\mathcal{C}$ be positive semi-definite with generating vector $c\in \mathbb{R}^n$. Suppose $x\in \mathbb{R}^n$ is the eigenvector of $\mathcal{C}$ corresponding to $\rho(\mathcal{C})$. Assume
\begin{equation}\label{e22} r_{\bar{i}}=\max_{1\leq i \leq n}r_i,~r_{\underline{i}}=\min_{1\leq i\leq n}r_i.
\end{equation}
Then, $R=r_{\bar{i}},~r=r_{\underline{i}}.$
\et
\proof Since Cauchy tensor $\mathcal{C}$ is positive semi-definite, from Theorem \ref{thm21}, all elements of $\mathcal{C}$ and $c$
are positive. By Definition \ref{def22}, we know that $\mathcal{C}$ is irreducible. Thus $x>0$ from Theorem 1.4 of \cite{Chang08}.
Without loss of generality, suppose
$$R=r_l,~r=r_s.$$
By the analysis before this theorem, it holds that
$$\bar{a}=\max_{1\leq i\leq n}c_i=c_s,~\underline{a}=\min_{1\leq i \leq n}c_i=c_l.$$

On the other side, by (\ref{e21}) of Definition \ref{def21}, we have
$$\begin{array}{rl}
\rho(\mathcal{C})x_{\bar{i}}^{m-1}=&(\mathcal{C}x^{m-1})_{\bar{i}}\\
=&\sum_{i_2,\cdots,i_m \in [n]}c_{\bar{i}i_2i_3\cdots i_m}x_{i_2} \cdots x_{i_m} \\
=&\sum_{i_2,\cdots,i_m \in [n]}\frac{x_{i_2} \cdots x_{i_m}}{c_{\bar{i}}+c_{i_2}+\cdots+c_{i_m}} \\
\leq& \sum_{i_2,\cdots,i_m \in [n]}\frac{x_{i_2} \cdots x_{i_m}}{\underline{a}+c_{i_2}+\cdots+c_{i_m}}\\
=&(\mathcal{C}x^{m-1})_l\\
=&\rho(\mathcal{C})x^{m-1}_l,\\
\end{array}$$
which implies that
$$x_{\bar{i}}=x_l$$
So we can take $\bar{i}=l$ and $R=r_{\bar{i}}$ holds.
Similarly, by (\ref{e21}) and (\ref{e22}), one has
$$\begin{array}{rl}
\rho(\mathcal{C})x_{\underline{i}}^{m-1}=&(\mathcal{C}x^{m-1})_{\underline{i}}\\
=&\sum_{i_2,\cdots,i_m \in [n]}c_{\underline{i}i_2i_3\cdots i_m}x_{i_2} \cdots x_{i_m} \\
=&\sum_{i_2,\cdots,i_m \in [n]}\frac{x_{i_2} \cdots x_{i_m}}{c_{\underline{i}}+c_{i_2}+\cdots+c_{i_m}} \\
\geq& \sum_{i_2,\cdots,i_m \in [n]}\frac{x_{i_2} \cdots x_{i_m}}{\bar{a}+c_{i_2}+\cdots+c_{i_m}}\\
=&(\mathcal{C}x^{m-1})_s\\
=&\rho(\mathcal{C})x^{m-1}_s,\\
\end{array}$$
which means that $x_{\underline{i}}=x_s$ and we can take $\underline{i}=s$. Thus $r=r_{\underline{i}}$ and the desired results follows.
\qed

\bt\label{thm26} Suppose even order Cauchy tensor $\mathcal{C}$ has positive generating vector $c\in \mathbb{R}^n$. Then $\mathcal{C}$ is positive definite
if and only if $r_1,r_2,\cdots,r_n$ are mutually distinct.
\et
\proof By conditions, all elements of $c$ are positive, so it is obvious that $r_1,r_2,\cdots,r_n$ are mutually distinct if and only if
$c_1,c_2,\cdots,c_n$ are mutually distinct. By Theorem \ref{thm23}, the desired conclusion follows.
\qed

Suppose $\mathcal{A}=(a_{i_1i_2\cdots i_m})$ and $\mathcal{B}=(b_{i_1i_2\cdots i_m})$ are two order $m$ dimension $n$ tensors, the Hadamard product
of $\mathcal{A}$ and $\mathcal{B}$ is defined as
$$\mathcal{A}\circ \mathcal{B}=(a_{i_1i_2\cdots i_m}b_{i_1i_2\cdots i_m}),$$
which is still an order $m$ dimension $n$ tensor.
\bt\label{thm27} Let $\mathcal{C}$ and $\mathcal{C}'$ be two order $m$ dimension $n$ Cauchy tensors with generating vectors $c=(c_1,c_2,\cdots,c_n)$ and $c'=(c_1',c_2',\cdots,c_n')$
respectively. Suppose $m$ is even. If Cauchy tensors $\mathcal{C}$ and $\mathcal{C}'$ are positive semi-definite, then $\mathcal{C}\circ \mathcal{C}'$ is a positive semi-definite tensor.
\et
\proof By conditions, from Theorem \ref{thm21}, it follows that
$$c>0,~c'>0.$$
Let $a'=c_1'+c_2'+\cdots+c_n'>0$.
For any vector $x\in \mathbb{R}^n$, by the definition of Hadamard product of tensors, we have
$$\begin{array}{rl} (\mathcal{C}\circ \mathcal{C}')x^m=& \sum_{i_1,\cdots,i_m \in [n]}c_{i_1i_2\cdots i_m}c_{i_1i_2\cdots i_m}'x_{i_1}x_{i_2} \cdots x_{i_m}\\
=&\sum_{i_1,\cdots,i_m \in [n]} \frac{x_{i_1}x_{i_2} \cdots x_{i_m}}{(c_{i_1}+c_{i_2}+\cdots+c_{i_m})(c_{i_1}'+c_{i_2}'+\cdots+c_{i_m}')} \\
=&\sum_{i_1,\cdots,i_m \in [n]}\displaystyle \int^1_0 t^{(c_{i_1}+c_{i_2}+\cdots+c_{i_m})(c_{i_1}'+c_{i_2}'+\cdots+c_{i_m}')-1}x_{i_1}x_{i_2} \cdots x_{i_m}dt\\
=&\sum_{i_1,\cdots,i_m \in [n]}\displaystyle \int^1_0 t^{c_{i_1}a'+c_{i_2}a'+\cdots+c_{i_m}a'-1}x_{i_1}x_{i_2} \cdots x_{i_m}dt\\
=&\displaystyle \int^1_0(\sum_{i\in [n]}t^{c_ia'-\frac{1}{m}}x_i)^m dt\\
\geq &0.
\end{array}$$
From the arbitrariness of vector $x$, we know that $\mathcal{C}\circ \mathcal{C}'$ is positive semi-definite.
\qed
\bc\label{corol24} Suppose $\mathcal{C}_1,~\mathcal{C}_2,\cdots,\mathcal{C}_l$ are order $m$ dimension $n$ positive semi-definite Cauchy tensors.
Assume $m$ is even. Then, $\mathcal{C}_1\circ\mathcal{C}_2\circ\cdots\circ\mathcal{C}_l$ is a positive semi-definite tensor.
\ec

By Theorems \ref{thm22} and \ref{thm23}, we have the following conclusion. We omit its proof since it is similar to the proof of
Theorem \ref{thm22}.
\bt\label{thm28} Let $\mathcal{C}$, $\mathcal{C}'$, $c$ and $c'$ be defined as in Theorem \ref{thm27}.
 If Cauchy tensors $\mathcal{C}$ and $\mathcal{C}'$ are positive definite, then $\mathcal{C}\circ \mathcal{C}'$ is a positive definite tensor.
\et

\setcounter{equation}{0}
\section{Inequalities for Cauchy Tensors}
In this section, we give several inequalities about the largest and the smallest H-eigenvalues of Cauchy tensors. The bounds for the largest
H-eigenvalues are given for positive semi-definite Cauchy tensors. Moreover, properties of Z-eigenvalues and Z-eigenvectors of odd order Cauchy tensors
are also shown.

It should be noted that a real symmetric tensor always has  Z-eigenvalues and an even order real symmetric tensor always has H-eigenvalues \cite{Qi05}. We denote the largest and smallest H-eigenvalues
of Cauchy tensor $\mathcal{C}$ by $\lambda_{max}$ and $\lambda_{min}$ respectively. When $\mathcal{C}$ is a positive semi-definite Cauchy tensor, then by the Perron-Frobenius theory of nonnegative tensors \cite{Chang08}, we have  $$\lambda_{max}=\rho(\mathcal{C}).$$
\bl\label{lema31} Assume $\mathcal{C}$ is a Cauchy tensor with generating vector $c$. If the entries of $c=(c_1,c_2,\cdots,c_n)$ have different signs, then,
$$\lambda_{min}\leq \frac{1}{m \max \{c_i|~c_i<0, i\in [n]\}}<0<\frac{1}{m \min \{c_i|~c_i>0, i\in [n]\}}\leq \lambda_{max}.$$
\el
\proof From Theorem 5 of \cite{Qi05}, we have
$$\lambda_{max}=\max \left\{\mathcal{C}x^m|~\sum_{i \in [n]}x_i^m=1,~x\in \mathbb{R}^n\right\}$$
and
$$\lambda_{min}=\min \left\{\mathcal{C}x^m|~\sum_{i \in [n]}x_i^m=1,~x\in \mathbb{R}^n\right\}.$$
Combining this with the fact that
$$\mathcal{C}e_i^m=\frac{1}{mc_i},~i\in [n],$$
we have the conclusion of the lemma.
\qed

Let $r$, $R$, $\bar{a}$ and $\underline{a}$ be defined as in Section 2. We have the following result.
\bt\label{thm31} Assume even order Cauchy tensor $\mathcal{C}$ has generating vector $c=(c_1,c_2,\cdots,c_n)$. Suppose $c>0$ and at least two elements of $c$ are
different.   Then
$$r+\frac{1}{m\bar{a}}\left(\sqrt{\frac{R}{r}}-1\right)< \lambda_{max} < R-\frac{1}{m\bar{a}}\left(1-\sqrt{\frac{r}{R}}\right).$$
\et
\proof Suppose $x\in \mathbb{R}^n$ is the eigenvector of $\mathcal{C}$ corresponding to $\lambda_{max}$. By conditions, Cauchy tensor $\mathcal{C}$ is an irreducible nonnegative tensor and it follows $x>0$ from Theorem 1.4 of \cite{Chang08}. Without loss of generality, let $x=(x_1,x_2,\cdots,x_n)$ and suppose $0<x_i\leq1,~i\in [n]$
such that
\begin{equation}\label{e31} x_s=\min_{i\in [n]}x_i>0,~x_l=\max_{j\in [n]}x_j=1.
\end{equation}
By Theorem \ref{thm25}, we have
$$R=r_l,~r=r_s$$
and $R> r$ since at least two entries of $c$ are not equal.

On the other side, by the definition of eigenvalues, from (\ref{e31}), one has
\begin{equation}\label{e32}
\begin{array}{rl} \lambda_{max}x_s^{m-1}=&(\mathcal{C}x^{m-1})_s\\
=&\sum_{i_2,\cdots,i_m \in [n]}c_{si_2\cdots i_m}x_{i_2}x_{i_3}\cdots x_m\\
\leq &\sum_{i_2,\cdots,i_m \in [n]}c_{si_2\cdots i_m} \\
=&r,
\end{array}
\end{equation}
and
\begin{equation}\label{e33}
\begin{array}{rl}\lambda_{max}=& \lambda_{max}x_l^{m-1} \\
=&(\mathcal{C}x^{m-1})_l\\
=&\sum_{i_2,\cdots,i_m \in [n]}c_{li_2\cdots i_m}x_{i_2}x_{i_3}\cdots x_m\\
\geq& x_s^{m-1}\sum_{i_2,\cdots,i_m \in [n]}c_{li_2\cdots i_m} \\
=&Rx_s^{m-1}.
\end{array}
\end{equation}
Thus, by (\ref{e32}) and (\ref{e33}), we have
$$
0<x_s^{m-1}\leq\frac{\lambda_{max}}{R}\leq\frac{r}{x_s^{m-1}R},
$$
which can be written as
$$
0<x_s^{m-1}\leq \sqrt{\frac{r}{R}}.
$$
Combining this with (\ref{e31}), we obtain
$$
\begin{array}{rl}\lambda_{max}=& \lambda_{max}x_l^{m-1} \\
=&(\mathcal{C}x^{m-1})_l\\
=&\sum_{i_2,\cdots,i_m \in [n], (l,i_2,\cdots,i_m)\neq(l,s,\cdots, s)}c_{li_2\cdots i_m}x_{i_2}x_{i_3}\cdots x_m +c_{ls\cdots s}x_s^{m-1}\\
<& \sum_{i_2,\cdots,i_m \in [n]}c_{li_2\cdots i_m}-c_{ls\cdots s}+c_{ls\cdots s}\sqrt{\frac{r}{R}}\\
<&R-\frac{1}{m\bar{a}}(1-\sqrt{\frac{r}{R}}),
\end{array}
$$
and
$$
\begin{array}{rl}\lambda_{max}=&\frac{1}{x_s^{m-1}}\sum_{i_2,\cdots,i_m \in [n]}c_{si_2\cdots i_m}x_{i_2}x_{i_3}\cdots x_m\\
=&\frac{c_{sl\cdots l}x_l^{m-1}}{x_s^{m-1}}+\frac{1}{x_s^{m-1}}\sum_{i_2,\cdots,i_m \in [n],(i_2i_3\cdots i_m)\neq (ll\cdots l)}c_{si_2\cdots i_m}x_{i_2}x_{i_3}\cdots x_m\\
=&\frac{c_{sl\cdots l}}{x_s^{m-1}}+\frac{1}{x_s^{m-1}}\sum_{i_2,\cdots,i_m \in [n],(i_2i_3\cdots i_m)\neq (ll\cdots l)}c_{si_2\cdots i_m}x_{i_2}x_{i_3}\cdots x_m\\
>&\sqrt{\frac{R}{r}}c_{sll\cdots l}+r-c_{sll\cdots l}\\
>&r+\frac{1}{m\bar{a}}(\sqrt{\frac{R}{r}}-1),
\end{array}
$$
from which we get the desired inequalities. \qed

In \cite{Qi05}, Qi called a real number $\lambda$ and a real vector $x\in \mathbb{R}^n$ a Z-eigenvalue of tensor $\mathcal{A}$ and a Z-eigenvector
of $\mathcal{A}$ corresponding to $\lambda$, it they are solutions of the following system:
\begin{equation}\label{e34}
\begin{cases}
\mathcal{A}x^{m-1}=\lambda x\\
x^Tx=1
\end{cases}.
\end{equation}
Next, we will give several spectral properties for odd order Cauchy tensors.

\bt\label{thm32} Suppose order $m$ dimension $n$ Cauchy tensor $\mathcal{C}$ has generating vector $c$. Let
$m$ be odd and $c>0$.
Assume $\lambda\in \mathbb{R}$ is a Z-eigenvalue of  $\mathcal{C}$ with Z-eigenvector $x=(x_1,x_2,\cdots,x_n)\in \mathbb{R}^n$.
If Z-eigenvalue $\lambda>0$, then $x\geq0$;
if Z-eigenvalue $\lambda<0$, then $x\leq0$.
\et
\proof By the condition $c>0$, we know that all entries of Cauchy tensor $\mathcal{C}$ are positive.
By definitions of Z-eigenvalue and Z-eigenvector, for any $i\in [n]$, we have that
\begin{equation}\label{e35}
\begin{array}{rl} \lambda x_i=&(\mathcal{C}x^{m-1})_i\\
=&\sum_{i_2,i_3,\cdots,i_m \in [n]}c_{ii_2\cdots i_m}x_{i_2}x_{i_3}\cdots x_m\\
=&\sum_{i_2,i_3,\cdots,i_m \in [n]}\frac{x_{i_2}x_{i_3}\cdots x_m}{c_i+c_{i_2}+\cdots+c_{i_m}}\\
=&\sum_{i_2,i_3,\cdots,i_m \in [n]}\displaystyle \int^1_0 t^{c_i+c_{i_2}+\cdots+c_{i_m}-1}x_{i_2}x_{i_3}\cdots x_m dt\\
=&\displaystyle \int^1_0 t^{c_i-\frac{1}{m}} (\sum_{j\in [n]}t^{c_j-\frac{1}{m}}x_j)^{m-1}dt.
\end{array}
\end{equation}
Since $m$ is odd, by (\ref{e35}), one has
$$\lambda x_i\geq 0,~~\mbox{for}~~ i\in [n],$$
which implies that $x\geq0$ when $\lambda>0$ and $x\leq o$ when $\lambda<0$.
\qed

\bt\label{thm33} Suppose Cauchy tensor $\mathcal{C}$ and its generating vector $c$ are defined as in Theorem \ref{thm32}.
If all entries of $c$ are mutually distinct, then $\mathcal{C}$ has no zero Z-eigenvalue.
\et
\proof By conditions, since entries of generating vector $c$ are mutually distinct, without loss of generality, suppose
$$0<c_1<c_2<\cdots<c_n.$$
We prove the result by contradiction. Suppose $\mathcal{C}$ has Z-eigenvalue $\lambda=0$ with Z-eigenvector $x\in \mathbb{R}^n$.
Then, by (\ref{e35}), for any $i\in [n]$, we have
$$\displaystyle \int^1_0 t^{c_i-\frac{1}{m}} \left(\sum_{j\in [n]}t^{c_j-\frac{1}{m}}x_j\right)^{m-1}dt\equiv0.$$
From properties of integration, one has
$$t^{c_i-\frac{1}{m}} \left(\sum_{j\in [n]}t^{c_j-\frac{1}{m}}x_j\right)^{m-1}\equiv0,~ t\in [0,1],$$
i.e.,
\begin{equation}\label{e36}
t^{c_i-\frac{1}{m}}\left(t^{c_1-\frac{1}{m}}x_1+t^{c_2-\frac{1}{m}}x_2+\cdots+t^{c_n-\frac{1}{m}}x_n\right)\equiv0,~t\in [0,1].
\end{equation}
By (\ref{e36}), we obtain
$$t^{c_1-\frac{1}{m}}x_1+t^{c_2-\frac{1}{m}}x_2+\cdots+t^{c_n-\frac{1}{m}}x_n\equiv0,~t\in (0,1],$$
which implies that
\begin{equation}\label{e37}x_1+t^{c_2-c_1}x_2+\cdots+t^{c_n-c_1}x_n\equiv0,~t\in (0,1].\end{equation}
Since $c_1,c_2,\cdots,c_m$ are mutually distinct, by the continuity property of operators and (\ref{e37}), it follows that
$$x_1=0.$$
Thus, the equation (\ref{e37}) can be written as
$$t^{c_2-c_1}x_2+\cdots+t^{c_n-c_1}x_n\equiv0,~t\in (0,1],$$
which is equivalent to
$$x_2+t^{c_3-c_2}x_3+\cdots+t^{c_n-c_2}x_n\equiv0,~t\in (0,1].$$
By the continuity property, we have
$x_2=0.$
Repeating the process above, we get
$$x_1=x_2=\cdots=x_n=0,$$
which is contradicting with the fact that $x$ is a Z-eigenvector corresponding to $\lambda=0$.   The desired conclusion follows.
\qed

\section{Final Remarks}
In this paper, we give several necessary and sufficient conditions for an even order Cauchy tensor to be positive semi-definite.
Some properties of positive semi-definite Cauchy tensors are presented. Furthermore, inequalities about the largest H-eigenvalue and
the smallest H-eigenvalue of Cauchy tensors are shown. At last, some spectral properties on Z-eigenvalues of odd order Cauchy tensors
are shown.

However, there are still some questions that we are not sure now. The Cauchy matrix can be combined with many other structured matrices to
form new structured matrices such as Cauchy-Toeplitz matrix and Cauchy-Hankel matrix \cite{Solak03,Tyr91,Tyr92}. Can we get the
type of Cauchy-Toeplitz tensors and Cauchy-Hankel tensors? If so, how about their spectral properties? What are the necessary and sufficient
conditions for their positive semi-definiteness?

\end{document}